\documentclass[12pt,a4paper,oneside]{amsart}
\usepackage{amsmath,amstext,amssymb,amsopn,amsthm,color}

\usepackage[T1]{fontenc} 
\theoremstyle{plain}
\newtheorem{thm}{Theorem}[section]
\newtheorem{theorem}{Theorem}[section]

\newtheorem{lemma}[thm]{Lemma}

\newtheorem{corollary}[thm]{Corollary}
\theoremstyle{definition}

\theoremstyle{remark}
\newtheorem*{rem*}{Remark}

\newtheorem{remark}{Remark}

\newcommand{\R}{\mathbb{R}}

\newcommand{\E}{\mathbb{E}}
\newcommand{\ex}{\mathbb{E}}
\newcommand{\p}{\mathbb{P}}
\newcommand{\pr}{\mathbb{P}}

\renewcommand{\H}{\mathbb{H}}

\renewcommand{\leq}{\leqslant}
\renewcommand{\geq}{\geqslant}

\renewcommand{\leq}{\leqslant}
\renewcommand{\geq}{\geqslant}

\def\({\left(}
\def\){\right)}
\def\[{\left[}
\def\]{\right]}
\def\<{\langle}
\def\>{\rangle}

\begin{document}

\title{On potential theory of hyperbolic Brownian motion with drift}
\thanks{The author was supported by Faculty of Pure and Applied Mathematics, Wroc{\l}aw University of Science and Technology}
\author{ Grzegorz Serafin}
\address{Grzegorz Serafin \\ Institute of Mathematics and Computer Science \\ Wroc{\l}aw University of Technology \\ ul. Wybrze{\.z}e Wyspia{\'n}\-skiego 27 \\ 50-370 Wroc{\l}aw, Poland}
\email{grzegorz.serafin@pwr.edu.pl}

\begin{abstract}
Consider the $\lambda$-Green function and the $\lambda$-Poisson kernel of 
a Lipschitz domain $U\subset \mathbb H^n=\left\{x\in\mathbb R^n:x_n>0\right\}$ for hyperbolic Brownian motion with drift. We provide several relationships that facilitate studying those objects and explain somehow theirs nature.  As an application, we yield  uniform estimates in case of sets of the form $S_{a,b}=\{x\in\H^n:x_n>a,x_1\in(0,b)\}$, $a,b>0$, which covers and extends existing  results of that kind.
\end{abstract}

\maketitle

\section{Introduction}
Hyperbolic Brownian motion (HBM) is a canonical  diffusion in the  real hyperbolic space with half of Laplace-Beltrami operator as its generator.  The process  is a natural counterpart of the classical Brownian motion  and plays a crucial role in probabilistic approach to the potential theory on hyperbolic space. On the other hand,  HBM is closely related to geometric Brownian motion and Bessel process \cite{BTFY}, \cite{Y2}. It has also some  applications to Physics \cite{GS}  and risk theory in Financial Mathematics  \cite{D}, \cite{Y1}.
Properties of HBM has been significantly developed in papers  \cite{BTF}, \cite{BTFY}, \cite{G}, \cite{M} and more. One of the main objects, in the context of potential theory on hyperbolic spaces, are  the $\lambda$-Green function and the $\lambda$-Poisson kernel of subdomains. They were recently intensively  studied for particular sets, see e.g., \cite{BGS,ByM1,BMZ,MS,S}. The main results of the paper (see Theorem \ref{thm:allis0}) implies that studying aforementioned objects  leads to the HBM with drift. For this reason, our approach is based on the process with drift from the very beginning.

We denote by $X^{(\mu)}=\{X^{(\mu)}(t)\}_{t\geq0}, \mu\in\R$,  the HBM with drift  on the half-space model $\H^n=\{x\in\R^n:x_n>0\}$ of the $n$-dimensional real hyperbolic space. The generator of the process is 
$\frac12\Delta_\mu$, where
\begin{equation}\label{eq:L-B}
\Delta_\mu = x_n^2\sum_{k=1}^n\frac{\partial^2 }{\partial x_k^2}-(2\mu-1)x_n\frac{\partial
}{\partial x_n}.
\end{equation}
Note  that $\Delta_{(n-1)/2}$ is the Laplace-Beltrami operator and $\mu=(n-1)/2$ corresponds therefore to the  standard HBM.  In the paper, we focus mostly on $\mu>0$, since the main motivation of studying HBM with drift, mentioned in the first paragraph, is related to positive values of $\mu$.  Furthermore, potential theories for opposite indices are associated to each other (see   Remark \ref{rem1} after Theorem \ref{thm:allis0}), which allows us to study only the positive ones.  

 Let us denote by $\tau_U^\mu=\inf\{t:X^{(\mu)}(t)\notin U\}$ the first exit time of the process from the domain $U$. The objective of the paper is to  examine the  $\lambda$-Green  function $G^{(\mu),\lambda}_U(x,y)$ and the $\lambda$-Poisson kernel $P^{(\mu),\lambda}_U(x,y)$ of  $U$, which are defined as follows
\begin{align}\label{defG}
G_U^{(\mu),\lambda}(x,y)=&\int_0^\infty e^{-\lambda t}\E^x\[t<\tau_U^\mu;X^{(\mu)}(t)\in dy\]dt/dy,\quad x,y\in U,\\\label{defP}
P_U^{(\mu),\lambda}(x,y)=&\,\E^x\[e^{-\lambda\tau_U^\mu};X^{(\mu)}\(\tau_U^\mu\)\in dy\]/dy, \quad x\in U, y\in \partial U.\end{align}
In the formula for $P_U^{(\mu),\lambda}(x,y)$ we assume additionally that $\tau_U^\mu<\infty$ a.s., if $\lambda>0$. Furthermore,
for $\lambda=0$ above objects became the Green function  and the Poisson kernel, which we  denote by  $G_U^{(\mu),0}(x,y)=G_U^{(\mu)}(x,y)$ and $P_U^{(\mu),0}(x,y)=P_U^{(\mu)}(x,y)$, respectively.   Those functions are fundamental objects in potential theory on $\H^n$. Precisely, they describe solutions for Dirichlet problem involving the operator $\Delta_\mu$. In particular, the Green function appears to be the inverse operator to $\Delta_\mu$.  The $\lambda$-Green function and $\lambda$-Poisson kernel take over the leading role when the operator $\Delta_\mu-\lambda I$ is considered. The $\lambda$-Green function may be then understood as the resolvent kernel for the operator $\Delta_\mu$, and the $\lambda$-Poisson kernel recovers $\lambda$-harmonic (($\Delta_\mu-\lambda I$)-harmonic) functions from boundary conditions.  The following relationships are provided in Theorem \ref{thm:allis0}:
\begin{equation}\label{allis0}
 G_U^{(\mu),\lambda}(x,y)=\(\frac{x_n}{y_n}\)^{\eta-\mu}G_U^{(\eta)}(x,y),\ \ 
 P_U^{(\mu),\lambda}(x,y)=\(\frac{x_n}{y_n}\)^{\eta-\mu}P_U^{(\eta)}(x,y),
\end{equation} 
where $\mu\in\R$ and $\eta=\sqrt{\mu^2+2\lambda}$. The main consequence of this result is that research on $\lambda$-Green functions and
$\lambda$-Poisson kernels can be reduced only to the case $\lambda=0$. Furthermore, studying above-mentioned objects for standard hyperbolic Brownian motion induces naturally introduction of HBM drift, which is a substantial motivation to study that process.

If  $\lambda>0$ and $\tau_U^\mu=\infty$ with positive probability,  the $\lambda$-Poisson kernel defined in a classical way  becomes degenerate. Indeed, since  $\{\tau_U^\mu=\infty\} =\{X_n\(\tau_U^\mu\)=0\}$ a.s., the right-hand side of (\ref{defP}) vanishes on the set $\partial U\cap P$, where $P=\{x\in\R^n:x_n=0\}$ and $\partial U$ is the boundary  in Euclidean metric (in $\R^n$) of $U$.  This effect is due to a specific behavior of $\lambda$-harmonic functions in a neighborhood of the set $P$. The  definition (\ref{defP}) does not take into consideration this behavior. To discuss this issue more precisely we recall an analytical interpretation of the $\lambda$-Poisson kernel as  an integral kernel  solving the Dirichlet problem.    Then we reformulate the problem and solve it by an integral kernel of the form corresponding to (\ref{allis0}).   

An another important result of the paper is Theorem \ref{XtoY} where we 
show that the Green function and the Poisson kernel  for HBM with drift  can be easily expressed by analogous objects  for Brown-Bessel diffusion. This  general method was introduced by Molchanov and Ostrovski \cite{MO}, see also \cite{O}. Finally, in Theorem \ref{thm:HtoE} we relate potential theory on $\H^n$ to the classical one on the Euclidean space $\R^{2n}$. However, that result concerns  HBM without drift and sets are being modified.   As  an application of general results we provide uniform  estimates of  the Green function and the Poisson kernel of the set $S_{a,b}=\{x\in\H^n:x_n>a,x_1\in(0,b)\}$, $a,b>0$. This set may seem very special, but studying it is motivated by geometry: hyperplanes $x_n=a$ are horocycles and hyperplanes $x_1=b$ are geodesics  in space $\H^n$. Moreover, constants in the estimates depend only on the dimension and parameter $\mu$, and manipulation of parameters $a$ and $b$ let us therefore recover and improve existing results for  sets such as $D_a=\{x\in\H^n:x_n>a\}$, $H=\{x\in\H^n:x_1>0\}$, $S_b=\{x\in\H^n:x_1\in(0,b)\}$ (see  \cite{BBM, BMR, MS, S}).

The paper is organized as follows. Preliminaries start with  a short  description of Bessel process and related objects. The process killed when exiting half-line $(a,\infty)$, $a>0$, is also considered. Afterwards  the hyperbolic space $\H^n$ and the  HBM with drift are introduced. In Section 3 we collect several relationships which simplify research on the $\lambda$-Green function and the $\lambda$-Poisson kernel of subdomains of $\H^n$. Section 4 is devoted to  estimates of the Green function and the Poisson kernel  of the set $S_{a,b}$. In Appendix, one can find an integral lemma which is intensively exploited in Section 4.

\section{Preliminaries}
\subsection{Notation}
 We present estimates using the following notation: for two positive  functions $f,g:\textit{X}\rightarrow(0,\infty)$ we write   $f \approx g$, if  there exists a constant $c>1$ such that \mbox{$1/c\leq f/g\leq c$} for every $x\in \textit{X}$. If the constant $c$ depends on an additional parameter, we write this parameter over the sing $\approx$. 
\subsection{Bessel process}
We denote by $R^{(\nu)}=\{R^{(\nu)}(t)\}_{t\geq 0}$ the  Bessel process with index $\nu<0$
starting from  $R^{(\nu)}(0)=x>0$. Nonnegative indices are also considered in literature, however, they are irrelevant from our point of view.  For $\nu\leq-1$ the point 0 is
killing and the process hits it a.s.. In the case $-1 < \nu< 0$ we
impose killing condition at 0. The transition density function of the process is given by (see \cite{BS} p.134)
\begin{equation}\label{eq:besseldens}
g^{(\nu)}(t,x,y)=\frac
y{t}\left(\frac{y}{x}\right)^\nu\exp\left(-\frac{x^2+y^2}{2t}\right)I_{|\nu|}\left(\frac{xy}{t}\right),\quad
\nu<0,\ x,y>0\/,
\end{equation}
where $I_{\nu}(z)$ is the modified Bessel function of the first kind.

Let  $B = \{B(t)\}_{t \geq 0}$ be the one-dimensional Brownian motion
starting from zero. Bessel process is related to the geometric Brownian
motion $\{x\exp(B(t)+\nu t)\}_{t\geq 0}$, $x>0$, by the Lamperti relation, which states
\begin{equation}\label{Lamperti}
\{x\exp\left(B(t)+\nu t\right)\}_{t\geq 0} \stackrel{d}{=}\left\{R^{(\nu)}\(A_x^{(\nu)}(t)\)\right\}_{t\geq0},
\end{equation}
where the integral functional $A_x^{(\nu)}(t)$ is defined by 
\begin{equation}\label{defA}
   A^{(\nu)}_{x}(t)=x^2\int_0^t\exp\left(2B(s)+2\nu s\right)ds.
\end{equation}
The density function $f_{x,t}^{(\nu)}(u,v)$ of a vector $\left(A^{(\nu)}_{x}(t),x \exp
\left(B(t)+\nu t\right)\right)$  was computed in \cite{Y1} and is given by
\begin{equation}\label{eq:jointdens}
f_{x,t}^{(\nu)}(u,v)=\left(\frac{v}{x}\right)^{\nu}e^{-\nu^2t/2}\frac1{uv}
\exp\left(-\frac{x^2+v^2}{2u}\right)\theta_{xv/u}(t),\ \ \ x,u,v,t>0.
\end{equation}
Here, the  function  $\theta_r(t)$ satisfies  (see \cite{Y3})
\begin{equation}\label{eq:Ltheta}
\int_0^\infty e^{-\lambda t}\theta_r(t)dt=I_{\sqrt{2\lambda}}(r).
\end{equation}
Note that the function $f_{x,t}^{(\nu)}(u,v)$  is also closely related to Hartman-Watson law (see \cite{HW}).

Bessel process with a negative index $\nu$ and starting from $x>a$, $a>0$, leaves the half-line $(a,\infty)$ with probability one. The transition density function of  the process killed on exiting  $(a,\infty)$   has been  estimates in  \cite{BM}, i.e. it holds
\begin{align}\nonumber g_a^{(\nu)}(t;x,y)\approx &\left[1\wedge \frac{(x-a)(y-a)}{t}\right]\left(1\wedge \frac{xy}{t}\right)^{|\nu|-\frac{1}{2}} \left(\frac{y}{x}\right)^{\nu+\frac{1}{2}}\frac{1}{\sqrt{t}}\exp\left(-\frac{(x-y)^2}{2t}\right)\\\label{kBdens}
\approx&\,\frac{(x-a)(y-a)}{t+(x-a)(y-a)}\(\frac{x^2}{t+xy}\)^{|\nu|-\frac{1}{2}} \frac{1}{\sqrt{t}}\exp\left(-\frac{(x-y)^2}{2t}\right),\quad x,y>a, t>0.
\end{align}
Furthermore, let us denote by $q^{(\nu)}_a(t;x)$, $a>0$, the density function  of the first hitting time of the point $a$ by the Bessel process. Its estimates are given in \cite{BMR} 
\begin{equation}\label{kBtime}
q^{(\nu)}_a(t;x)\approx\frac{(x-a)}{t^{3/2}} \frac{ x^{2|\nu|-1} }{\(t+ ax\)^{|\nu|-1/2}}\exp\left(-\frac{(x-y)^2}{2t}\right)\/,\quad t>0,x>a,\nu< 0\/.
\end{equation}
In fact, authors of \cite{BMR} have  made a simple mistake in the  formulation of this result and the above-given correct formula  differs slightly from the original one. Precisely, Theorems $4$ and $8$ in \cite{BMR} providing estimates when $x>t$ and when $x<t$, respectively, are correct, but the formula $(15)$ combining those theorems  is wrong: there should be $(t+x)^{|\mu|-1/2}$ instead of $t^{|\mu|-1/2}+x^{|\mu|-1/2}$ in the denominator.

\subsection{Hyperbolic space and hyperbolic Brownian motion with drift}
We consider the half-space model of the real hyperbolic space
$$\H^n=\left\{x\in\R^n:x_n>0\right\}, \ \ n=1,2,3,...\, . $$
The  formula for hyperbolic distance is given by
\begin{align}
\label{coshdxy}
 \cosh d_{\H^n}(x,y)&=\(1+\frac{\left|x-y\right|^2}{2x_ny_n}\),\ \ \ \ \ x,y\in\H^n.
\end{align}
  The unique, up to a constant factor, second order elliptic differential
operator on $\H^n$, annihilating constant functions, which is invariant under
isometries of the space is the Laplace-Beltrami operator $\Delta_{(n-1)/2}$ (cf.(\ref{eq:L-B})).

Hyperbolic Brownian motion with drift is a process $X^{(\mu)} = \{X^{(\mu)}(t)\}_{t\geq0}$ starting from $ X^{(\mu)}(0)=x\in\H^n$ which generator is $\frac12\Delta_\mu$.  The parameter $\mu$ is called  an index and the drift is equal to $\mu-\frac{n-1}{2}$. Note that for $\mu=\frac{n-1}{2}$ we obtain standard HBM (without drift). 

Let us denote by $B(t)=\left(B_1(t),...,B_n(t)\right)$ the classical Brownian motion in $\R^n$ starting from $(x_1,...,x_{n-1},0)$. Then the HBM with drift may be represented in terms of the process $B(t)$ as follows (see \cite{BBM})
\begin{equation}\label{eq:characterisation}
X^{(\mu)}(t)\stackrel{d}{=}\left(B_1(A^{(-\mu)}_{x_n}(t)),...,B_{n-1}(A^{(-\mu)}_{x_n}(t)),x_n \exp
\left(B_n(t)-  \mu t\right)\right).
\end{equation}
Here, the integral functional $A^{(-\mu)}_{x_n}(t)$, defined by (\ref{defA}),  is associated with $B_n(t)$.  In addition, using \emph{Lamperti relation}, we get
\begin{eqnarray}
\label{eq:subordinated Y}
\left\{X^{(\mu)}(t);t\geq0\right\}\stackrel{(d)}{=}\left\{Y\left(A^{(-\mu)}_{x_n}(t)\right);t\geq0\right\},
\end{eqnarray}
where 
\begin{equation}
\label{eq:defY}Y(t)=\left(B_1\left(t\right),...,B_{n-1}\left(t\right),R^{(-\mu)}\left(t\right)\right),\end{equation}
and the process $R^{(-\mu)}(t)$ is the  Bessel process with index $-\mu$ starting from $x_n$ and independent of the process $\left(B_1(t),...,B_{n-1}(t)\right)$.

\section{General results}
\subsection{Reduction to $\lambda=0$}
In this subsection we provide  precise relationships which  bond  $\lambda$-Green function and $\lambda$-Poisson kernel with analogous objects for $\lambda=0$ and for the process with different drift. It lets us reduce $\lambda$-potential theory to the case $\lambda=0$. The only cost we pay is mentioned change of drift of the process. 
\begin{theorem}\label{thm:allis0}
Let $\mu\in \R$ and  $U$ be a  domain in $\H^n$ and $\lambda\geq0$. We have
\begin{equation}\label{lGsimplification}G^{(\mu),\lambda}_U(x,y)=\(\frac{x_n}{y_n}\)^{\eta-\mu}G^{(\eta)}_U(x,y),\ \ \ \ x,y\in U,
\end{equation}
where $\eta=\sqrt{\mu^2+2\lambda}$. If, additionally, $\tau_U^{\mu}<\infty$ a.s., we get
\begin{equation}\label{lPsimplification}
P^{(\mu),\lambda}_U(x,y)=\(\frac{x_n}{y_n}\)^{\eta-\mu}P^{(\eta)}_U(x,y),\ \ \ \ x\in U, y\in\partial U.
\end{equation}
\end{theorem}

\begin{proof}
The last coordinate of the process $X^{(\mu)}(t)$ can be expressed in the form $X_n^{(\mu)}(t) = x_n\exp\left(W^{(\eta-\mu)}(t)-\eta
t\right)$, where $W^{(\eta-\mu)}(t)=B_n(t)+(\eta-\mu)t$ and $B_n$ is a one-dimensional Brownian motion. By virtue of the Girsanov theorem,  the process $\{W^{(\eta-\mu)}(t)\}_{0\leq t\leq T}$  is , for every $T>0$,  a standard Brownian motion with respect to the measure  $Q_T$ given by
\begin{equation}\label{M(T)}
  \frac{dQ_T}{d\pr}=\exp\left((\mu-\eta)B_n(T)-\frac12(\eta-\mu)^2T\right)=M(T).
\end{equation}
It implies that the process $\{ X^{(\mu)}(t)\}_{ 0\leq t \leq T}$ considered with respect to the measure $Q_T$ is a hyperbolic Brownian motion with drift with index $\eta$. Hence, for every $t\leq T$ and every Borel set $A\subset U$ we get
$$\ex^x\[t<\tau^\eta_U;X^{(\eta)}(t)\in A\]=\ex^x\[t<\tau^\mu_U;M(T);X^{(\mu)}(t)\in A\].$$
Observe that $M(T)$ is a $\mathcal F_t$-martingale and it may be rewritten as
\begin{align*}
M(T)&=\[\exp\left(B_n\left(t\right)-\mu t\right)\]^{\mu-\eta}e^{-(\eta^2-\mu^2)t/2}=x_n^{\eta-\mu}e^{-\lambda t}\[X_n^{(\mu)}(t)\]^{\mu-\eta}.
\end{align*}
Furthermore, let us denote by $\mathcal F_t$ the $\sigma$-field generated by $\{X_s^{(\mu)}\}_{0\leq s\leq t}$. The set $\left\{t<\tau_U^\mu\right\}$ is then
$\mathcal F_t$-measurable and we get for $t\leq T$
\begin{align*}
\ex^x\[t<\tau^\eta_U;X^{(\eta)}(t)\in A\]&=\ex^x\[\ex^x\[t<\tau^\mu_U;M(T);X^{(\mu)}(t)\in A\left|\mathcal F_t\right.\]\]\\
&=\ex^x\[t<\tau^\mu_U;\ex^x\[M(T)\left|\mathcal F_t\right.\];X^{(\mu)}(t)\in A\]\\
&=\ex^x\[t<\tau^\mu_U;M(t);X^{(\mu)}(t)\in A\]\\
&=x_n^{\eta-\mu}e^{-\lambda
t}\,\ex^x\[t<\tau^\mu_U;\left(X_n^{(\mu)}(t)\right)^{\mu-\eta};X^{(\mu)}(t)\in A\].
\end{align*}
Since there is no upper bound of $T$, the above-given equalities work for every $t\geq0$ and consequently
\begin{align*}
\int_AG^{(\eta)}_U(x,y)dy=&\int_0^\infty \ex^x\[t<\tau^\eta_U;X^{(\eta)}(t)\in A\]dt\\
=&\,x_n^{\eta-\mu}\int_0^\infty e^{-\lambda
t}\,\ex^x\[t<\tau^\mu_U;\left(X_n^{(\mu)}(t)\right)^{\mu-\eta};X^{(\mu)}(t)\in A\]dt\\
=&\,x_n^{\eta-\mu}\int_Ay_n^{\mu-\eta}G^{(\mu),\lambda}_U(x,y)dy,
\end{align*}
which proves the formula (\ref{lGsimplification}). Let us now focus on the proof of the latter assertion of the theorem. Similarly as before, for every Borel $C\subset\partial U$  we get
\begin{align*}
\ex^x\left[\tau^\eta_U<T;X^{(\eta)}({\tau^\eta_U})\in C
\right]&=\ex^x\left[\tau^\mu_U<T;M(T);X^{(\mu)}({\tau^\mu_U})\in
C \right]\\
&=\ex^x\left[\tau^\mu_U<T;\ex^x\left[M(T)\left|\mathcal
F_{\tau^\mu_U}\right.\right];X^{(\mu)}\(\tau^\mu_U\)\in C\right],
\end{align*}
where 
$$\mathcal F_{\tau^\mu_U}=\left\{A\in\mathcal F_\infty:\forall(t\geq0)A\cap\{\tau^\mu_U<t\}\in\mathcal
F_t\right\}.$$
Using Doob's optional stopping theorem we obtain
\begin{align*}
\ex^x&\left[\tau^\eta_U<T;X^{(\eta)}({\tau^\eta_U})\in C \right]=\ex^x\left[\tau^\mu_U<T;M(\tau^\mu_U);X^{(\mu)}({\tau^\mu_U})\in
C\right]\\
&=x_n^{\eta-\mu}\ex^x\left[\tau^\mu_U<T;(X^{(\mu)}_n\left(\tau^\mu_U\right))^{\mu-\eta}e^{-\lambda\tau^\mu_U};X^{(\mu)}({\tau^\mu_U})\in
C\right].
\end{align*}
The next step is to take a limit as  $T\rightarrow\infty$. By the assumption  $\tau_U<\infty$ a.s.,  monotone convergence theorem gives us
\begin{equation*}
\ex^x\left[X^{(\eta)}({\tau^\eta_U})\in C \right]=x_n^{\eta-\mu}\ex^x\left[(X^{(\mu)}_n\left(\tau^\mu_U\right))^{\mu-\eta}e^{-\lambda\tau^\mu_U};X^{(\mu)}({\tau^\mu_U})\in C\right],
\end{equation*}
which is equivalent to  (\ref{lPsimplification}).
\end{proof}
\begin{remark}\label{rem1}
An another significance of the above theorem is that is shows some kind of duality for potential theories for opposite values of the parameter $\mu$. Taking $\lambda=0$ and $-\mu$ instead of $\mu$ we get for $\mu>0$
\begin{equation}\label{-mu}
 G_U^{(-\mu)}(x,y)=\(\frac{x_n}{y_n}\)^{2\mu}G_U^{(\mu)}(x,y),\ \ 
 P_U^{(-\mu)}(x,y)=\(\frac{x_n}{y_n}\)^{2\mu}P_U^{(\eta)}(x,y).
\end{equation} 
 
\end{remark} We  turn now to the case when $\p^x\(\tau_U^\mu=\infty\)>0$. Since $X^{(\mu)}_n\(\infty\)=0$ (cf. (\ref{eq:characterisation})), the right-hand side of (\ref {defP}) vanishes at  $y_n=0$ for $\lambda>0$. This situation is singular, especially from the analytical point of view. Namely, the $\lambda$-Poisson kernel is supposed to solve the Dirichlet problem with  a given boundary condition. In our situation,  the condition on the set $P_U^{(\mu),\lambda}(x,y)=0$ (=$\partial U\cap\{y_n=0\}$)  has no influence on  behaviour of the solution  in the neighbourhood   of  that set.  However, some examples show that this behaviour is relevant and setting boundary conditions only on the set  $\partial U\cap\{y_n>0\}$ results in infinite number of solutions. Furthermore, we can observe that multiplying the right-hand side of (\ref{lPsimplification}) by $y_n^{\eta-\mu}$ and enlarging the set $U$ we obtain a nontrivial object. Finally, if a function  $f(x)$ is $\lambda$-harmonic  for the operator $\frac12\Delta_\mu$ (i.e. $\frac12\Delta_\mu f(x)=\lambda f(x)$), then the function $x_n^{\eta-\mu}f(\mu)$, where $\eta=\sqrt{2\lambda+\mu^2}-\mu$, is harmonic for the operator $\frac12\Delta_\eta$, which comes from the following
\begin{align}\nonumber
\frac12\Delta_\eta \left(x_n^{\eta-\mu}f(x)\right)
=&\ x_n^{\eta-\mu}\frac12\Delta_\eta f(x)+(\eta-\mu)x_n^{\mu-\eta+1}\frac{\partial f}{\partial
x_n}(x)-\frac{\eta^2-\mu^2}{2}x_n^{\mu-\eta}f(x)\\\label{eq:deltadelta}
=&\ x_n^{\mu-\eta}\frac12\Delta_\mu f(x)-\lambda
x_n^{\mu-\eta}f(x)=0.
\end{align}

One can show (using e.g., Theorem 4 in \cite{T}), that every continuous and bounded function on a Lipschitz domain $U$, which is harmonic for  $\Delta_\eta$, $\eta>0$,   has a limit at the boundary of $U$.  All this leads us to the following modified Dirichlet problem:

Set a Lipchitz domain $U\subset\H^n$,   $f\in \mathcal{C}_b(\partial U)$ and $\lambda>0$. Find a function 
$u\in\mathcal{C}^2(U)$
satisfying differential equation
	\begin{equation}
\label{dirichlet1}
	\left(\frac12\Delta _{\mu}u\right)(x)=\lambda u(x),\quad  x \in U\/,
\end{equation}
 such that the function $x_n^{\sqrt{2\lambda+\mu^2}-\mu}u(x)$ is bounded and
		\begin{equation}\label{dirichlet2}
\lim_{\substack{x\rightarrow z\\ x\in U}}x_n^{\sqrt{2\lambda+\mu^2}-\mu}u(x)=f(z)\/,\quad z\in
\partial U\/.
\end{equation}

\begin{theorem}\label{thm:solutionofDP}
The function $u$ satisfying (\ref{dirichlet1}) and (\ref{dirichlet2})  is unique and given by

$$
u(x)=x_n^{\mu-\eta}\int_{\partial U}f(y)P^{(\eta)}_U(x,y)dy,
$$
 where $\eta = \sqrt{\mu^2+2\lambda}$.
\end{theorem}
\begin{remark}
According to this theorem, we can  treat the function $x_n^{\mu-\eta}P^{(\eta)}_U$ as a kind of $\lambda$-Poisson kernel.  It does not cover the formula  for the $\lambda$-Poisson kernel from theorem \ref{thm:allis0}, but the only  difference is  the factor $y_n^{\mu-\eta}$.
\end{remark}

\begin{proof}
Define a function
$h(x)= x_n^{\eta-\mu}u(x) = \ex^x\left[f\left(X^{(\eta)}({\tau^\eta_U})\right)\right]$. It is bounded by $\left\|f\right\|_\infty$ and, according to the stochastic continuity of the process $X^{(\mu)}$, satisfies  condition
(\ref{dirichlet2}). Since $P^{(\eta)}_U(x,y)$ is the standard Poisson kernel for the process $X^{(\eta)}(t)$, we have $\Delta_\eta h(x)=0$. Thus, similarly as in  (\ref{eq:deltadelta}), we get $\frac12\Delta_\mu u(x)=\lambda u(x)$.

To prove the uniqueness of the solution let us consider a sequence of bounded, in hyperbolic metric, sets such that $U_m\nearrow U$. For every  $m$ the function $u_{\upharpoonright U_m}$ 
satisfies (\ref{dirichlet1}) and (\ref{dirichlet2}) for  $U_m$ instead of $U$ and for $f=u_{\upharpoonright \partial
U_m}\in\mathcal{C}_b(\partial U_m)$ and it is the unique function of this property (see \cite{Fo}). Moreover, we have (see \cite{KS} prop. 7.2, p. 364)
\begin{equation*}
  u_{\upharpoonright U_m}(x) = \ex^x\left[e^{-\lambda{\tau}^\mu_{U_m}}u\left(X^{(\mu)}({\tau^\mu_{U_m}})\right)\right]\/.
\end{equation*}
Hence, by Theorem  \ref{thm:allis0}, we get
\begin{equation*}
u(x) = \ex^x\left[e^{-\lambda{\tau}^\mu_{U_m}}u\left(X^{(\mu)}({\tau_{U^\mu_m}})\right)\right]
=
x_n^{\mu-\eta}\ex^x\left[(X^{(\eta)}({\tau^\eta_{U_m}}))^{\eta-\mu}u\left(X^{(\eta)}({\tau^\eta_{U_m}})\right)\right], \ \ \ x\in U_m\/.
\end{equation*}
As $m$ tends to infinity, by the Lebesgue's dominated convergence theorem we obtain
$$
u(x)=x_n^{\mu-\eta}\ex^x\left[f\left(X^{(\eta)}({\tau^\eta_U})\right)\right].
$$
\end{proof}
\subsection{Representations involving other processes}
Let us define  the Green function $G^Y_U(x,y)$ and the Poisson kernel $P^Y_U(x,y)$ of the set $U\subset \H^n$ for the process Brown-Bessel diffusion  $Y(t)$ (see (\ref{eq:defY})) analogously as for the HBM with drift, i.e.
\begin{align*}
&P^Y_U(x,y){=}\p^x\(Y\left(\tau_U^Y\right)\in dy \),\ \ \ \ x\in U, y\in\partial U,\\
&G_U^Y(x,y)=\int_0^\infty \E^x\[t<\tau_U^Y;Y(t)\in dy\]dt,\ \ \ x,y\in U.
\end{align*}
The next lemma  let us study these objects instead of theirs counterparts    for HBM with drift. The main advantage of  this result comes from independence of coordinates of the process $Y(t)$ and from the fact that this process is relatively well known. 
\begin{theorem}\label{XtoY}
For any  domain $U\subset\H^n$  we have
\begin{align*}
&(i)\ P_U^{(\mu)}(x,y)=P_U^Y(x,y),\\
&(ii)\ G_U^{(\mu)}(x,y)=\frac1{y_n^2}G^Y_U(x,y).
\end{align*}
\end{theorem}
\begin{remark}
The first assertion may be find in \cite{MS}, however, the proof is short so we repeat it for convenience of the Reader. The other assertion is proved in case of the set $D_a=\{x\in \H^n:x_n>a\}$ in \cite{BBM}, but the below-given proof    is much simpler and covers general sets.
\end{remark}

\begin{proof}
According to the representation (\ref{eq:subordinated Y}), the process $Z(t)=Y\left(A^{(-\mu)}_{x_n}(t)\right)$ is a HBM with drift.
Since the functional $A^{(-\mu)}_{x_n}(t)$ is continuous and increasing a.s., we have
$\tau_U^Y=A^{(-\mu)}_{x_n}(\tau_U^\mu)$ a.s.. Thus
$$X^{(\mu)}\left(\tau^\mu_U\right)\stackrel{d}{=}Z\left(\tau_U^\mu\right)=Y\left(A^{(-\mu)}_{x_n}(\tau_U^\mu)\right)\stackrel{a.s.}{=}Y\left(\tau_U^Y\right).$$
Denote by $p^{(\mu)}(t;x,y)$ the transition density function (with respect to the Lebesgue measure) of the process $Z(t)$. By the Hunt formula and the  Fubini-Tonelli theorem we have
\begin{align}\nonumber
\int_0^\infty \E^x\[t<\tau_U^\mu;Z(t)\in dy\]dt=&\int_0^\infty p^{(\mu)}(t;x,y)-\E^x\[t>\tau_U^\mu;p^{(\mu)}\(t-\tau_U^\mu, Z(\tau_U^Z),y\)\]dt \\\nonumber
=&\int_0^\infty p^{(\mu)}(t;x,y)dt-\E^x\[\int_{\tau_U^Z}^\infty p^{(\mu)}\(t-\tau_U^Z, Z(\tau_U^\mu),y\)dt\]\\
\label{PotentialZY}
=&\int_0^\infty p^{(\mu)}(t;x,y)dt-\E^x\[\int_0^\infty p^{(\mu)}\(t, Z(\tau_U^\mu),y\)dt\].
\end{align}
Using representation  (\ref{eq:characterisation}) and formulae (\ref{eq:jointdens}), (\ref{eq:Ltheta}) we get
\begin{align*}
\int_0^\infty p^{(\mu)}(t;x,y)dt=&\int_0^\infty \int_0^\infty\frac{1}{(2\pi u)^{(n-1)/2}}e^{-(\tilde x-\tilde y)^2/2u}f^{(\mu)}_{x_n,t}(u,y_n)du\,dt\\
=&\int_0^\infty \frac{1}{(2\pi u)^{(n-1)/2}}e^{-(\tilde x-\tilde y)^2/2u}e^{-\frac{x_n^2+y_n^2}{2u}}I_{\mu}\(\frac{x_ny_n}u\)du\\
=&\frac1{y_n^2}\int_0^\infty \frac{1}{(2\pi u)^{(n-1)/2}}e^{-(\tilde x-\tilde y)^2/2u}g^{(-\mu)}(u;x_n,y_n)du,
\end{align*}
where $g^{(\nu)}(u;x,y)$ is the transition density function of a Bessel process with index $\nu$ starting from $x$. We identify the function under the last integral  as the transition density function of the process $Y(t)$. Since $ Z\left(\tau^\mu_U\right)\stackrel{d}{=}Y\left(\tau_U^Y\right)$ holds and the property (\ref{PotentialZY}) can be derived also for the process $Y(t)$,  we obtain the statement $(ii)$. 
\end{proof}
The above theorem, together with scaling properties of the standard Brownian motion and  Bessel process, gives us the  following scaling properties of the Green function and  the Poisson kernel for HBM with drift.

\begin{corollary}\label{scaling}
For any  domain $U\subset \H^n$ and $a>0$ we have
\begin{align}\label{skalowanieG}
G^{(\mu)}_{aU}(x,y)=&\frac1{a^n}G^{(\mu)}_U\(\frac xa,\frac ya\),\\\label{skalowanieP}
P^{(\mu)}_{aU}(x,y)=&\frac1{a^{n-1}}P^{(\mu)}_U\(\frac xa,\frac ya\).
\end{align}
\end{corollary}

The last theorem in this section exhibits that the Green function and the Poisson kernel of a set $U$ for the HBM (without drift) in $\H^n$ may be  derived from theirs counterparts for a somewhat modified set and for classical Brownian motion in $\R^{2n}$. This shows that studying potential theory for HBM may be reduced to the classical potential theory. Note that such results in case of a class of  tube domains were obtained in \cite{S}.

For $A\subset\H^n$ we define 
$$ A^+:=\{x\in\R^{2n}:\big(x_1,x_2,...,x_{n-1},|(x_n,x_{n+1},...,x_{2n})|\big)\in A\}\subset\R^{2n},$$   
and
$$ x^+=(x,0,...,0)\in\R^{2n}.$$
The form of the set $A^+$ may be, in general,  slightly discouraging, but in many cases it is not very complicated, e.g. for the set $\{x\in\H^n: x_n<a\}$, a>0, or for tube domains. Furthermore, we define by $G^B_U(x,y)$ and $P^B_U(x,y)$ the Green function  and the Poisson kernel, respectively, of the set $U\subset\R^{2n}$ for the classical Brownian motion in $\R^{2n}$.
\begin{theorem}\label{thm:HtoE}
For any open set  $A\in\H^n$ we have
$$G_A(x,y)=\frac{x_n^{n-1}}{y_n^{n+1}}\int_{|(w_1,...,w_{n+1})|=y_n}G^B_{A^+}\(x^+,(\tilde y,w)\)d\sigma(w),\ \ \ \ \ \ x,y\in A,$$
where $\tilde y=(y_1,...,y_{n-1})$. Furthermore, if $\partial A\cap (\R^{n-1}\times\{0\})=\phi$, then we have
$$P_A(x,y)=\(\frac{x_n}{y_n}\)^{n-1}\int_{|(w_1,...,w_{n+1})|=y_n}P^B_{A^+}\(x^+,(\tilde y,w)\)d\sigma(w),\ \ \ x\in A, y\in\partial A.$$
\end{theorem}
\begin{proof}

By the formula (\ref{-mu}) and Theorem \ref{XtoY} we get 
\begin{align*}
G^{\(-\frac{n-1}2\)}_U(x,y)&=\(\frac{x_n}{y_n}\)^{n-1}G^{\(\frac{n-1}2\)}_U(x,y)=\frac{x_n^{n-1}}{y_n^{n+1}}G^{Y}_U(x,y),\\
P^{\(-\frac{n-1}2\)}_U(x,y)&=\(\frac{x_n}{y_n}\)^{n-1}P^{\(\frac{n-1}2\)}_U(x,y)=\(\frac{x_n}{y_n}\)^{n-1}G^{Y}_U(x,y),
\end{align*}
where $Y(t)=\left(B_1\left(t\right),...,B_{n-1}\left(t\right),R^{\(\frac{n-1}2\)}\left(t\right)\right)$. Since the Bessel process $R^{\(\frac{n-1}2\)}$ may be interpreted as a norm of $n+1$-dimensional standard Brownian motion, we obtain
$$Y\stackrel{d}{=}\(B_1,B_2,...,B_{n-1},\|(B_n,...,B_{2n})\|\),$$
where $B=(B_1,...,B_{2n})$ is a $2n$-dimensional Brownian motion starting from $x^+$. Let $f$ be a positive function on $\H^n$. Following the convention that $G_A(x,\cdot)$ vanishes on $A^c$ we get
\begin{align*}
\int_{A}&G_A^Y(x,y)f(y)dy=\int_0^\infty \E^x\[f(Y);t<\tau_A^Y\]dt\\
&=\int_0^\infty \E^{x^+}\[f\(B_1(t),...,B_{n-1}(t),\|\(B_n(t),...,B_{2n}(t)\)\|\);t<\tau_{A^+}^B\]dt\\
&=\int_{\R^{2n}}f\(y_1,...,y_{n-1},\sqrt{z_1^2+...+z_{n+1}^2}\)G^B_{A^+}(x^+, (y_1,...,y_{n-1},w))dy_1...dy_{n-1}dw_1...dw_{n+1}\\
&=\int_{\H^n}f\(y_1,...,y_{n-1},y_n\)\int_{|(w_1,...,w_{n+1})|=y_n}G^B_{A^+}\(x^+,(\tilde y,w)\)d\sigma(w)dy_1...dy_{n},
\end{align*}
where $\sigma$ denotes the $n+1$-dimensional spherical measure. Furthermore, let $\tau_{A^+}^B$ be the first exit time of the Brownian motion $B(t)$ from the set $A^+$. Then for every   positive function $g$ on $\partial A$ we have
\begin{align*}
\E^x&\[g(Y(\tau_A^{Y}))\]= \E^{x^+}\[g\(B_1(\tau_{A^+}^B),...,B_{n-1}(\tau_{A^+}^B),\|\(B_n(\tau_{A^+}^B),...,B_{2n}(\tau_{A^+}^B)\)\|\)\]\\
&=\int_{\partial A^+}g\(y_1,...,y_{n-1},z_1,...,z_{n+1}\)P^B_{A^+}\(x^+,(\tilde y,z)\)dy_1...dy_{n-1}dz\\
&=\int_{\partial A}g\(y_1,...,y_{n-1},y_n\)\int_{|(w_1,...,w_{n+1})|=y_n}P^B_{A^+}\(x^+,(\tilde y,w)\)d\sigma(w)dy.
\end{align*}
Here, $dy$ stands for the induced Lebesgue measure on $\partial A$.
\end{proof}

\section{Estimates}

For $a,b>0$ we define
$$S_{a,b}=\{x\in\H^n:x_n>a, x_1\in(0,b)\}.$$
Studying this kind of sets is motivated by the hyperbolic geometry. The set $S_{a,b}$ is bounded by three hyperplanes: $P_1=\{x\in\H^n:x_1=0\}$, $P_2=\{x\in\H^n:x_1=b\}$ and $P_3=\{x\in\H^n:x_n=a\}$. Symmetries  with respect to hyperplanes $P_1$ and $P_2$ are isometries in $\H^n$; the set $P_3$ is a horocycle. In this section we estimate the Green function and the Poisson kernel of $S_{a,b}$ uniformly with respect to space variables as well as to  parameters $a$ and $b$. This lets us provide estimates for some other sets that may be obtained from $S_{a,b}$ by manipulation of values of the parameters .

By $\delta_u(w)=w\wedge (u-w)$, $u>0$, $w\in(0,u)$, we denote the Euclidean distance between $w$ and a compliment of the interval $(0,u)$. We clearly have $\delta_u(w)\,{\approx}\,w(u-w)/u$. 
Moreover, for $x\in\R^n$ and $a>0$ we define
\begin{equation}\label{x-a}x^{\downarrow a}=(x_1,...,x_{n-1},x_n-a).\end{equation}
\begin{theorem}\label{estG}
For $x,y\in S_{a,b}$ we have
$$G^{(\mu)}_{S_{a,b}}(x,y)\stackrel{\mu,n}{\approx}\frac{x_n^{\mu-1/2}}{y_n^{\mu+3/2}}\frac{e^{-\frac{\pi}b\left|x-y\right|}}{\left|x-y\right|^{n}}\frac{\[\delta_b(x_1)\delta_b(y_1)\]\wedge\left|x-y\right|^2}{\(\frac1b\left|x-y\right|+\cosh\rho_a\)}\frac{\(1+\frac1b\left|x-y\right|\)^{n/2+\mu+3/2}}{\(\frac1b\left|x-y\right|+\cosh\rho\)^{\mu-1/2}},$$
where $\rho_a$ is a hyperbolic distance between $x^{\downarrow a}$ and $y^{\downarrow a}$.
\end{theorem}
\begin{proof}Scaling property follows $G^{(\mu)}_{S_{a,b}}(x,y)=\frac1{b^n}G^{(\mu)}_{S_{a/b,1}}\(\frac xb,\frac yb\)$, hence it is enough to consider $b=1$. Furthermore, statement $(ii)$ in  Lemma \ref{XtoY} gives us
$$G_{S_{a,1}}^{(\mu)}(x,y)=\frac1{y_n^2}\int_0^\infty
j(t; x_1,y_1)\frac{\exp\(-\frac1{2t}\sum_{k=2}^{n-1}(x_k-y_k)^2\)}{(2\pi
t)^{(n-2)/2}}g^{(-\mu)}_a(t;x_n,y_n)dt,$$
where $g^{(-\mu)}_a(t;x_n,y_n)$ is the transition density of a Bessel process with index $-\mu$ killed on exiting $(a,\infty)$ and $j(t; x_1,y_1)$ is the transition density function of a one-dimensional Brownian motion killed on exiting the interval $(0,1)$. Estimates of the function $j(t; x_1,y_1)$ are given in Theorem 5.4 in \cite{PSZ} (cf. \cite{Fe}, (5.7) p. 341 and \cite{PSZ}, Theorem 2.2):
\begin{align}\nonumber
j(t; x_1,y_1)\approx &\(1\wedge\frac{x_1y_y}{t}\)\(1\wedge\frac{(1-x_1)(1-y_1)}{t}\)\frac{1+t^{5/2}}{\sqrt t}e^{-\pi^2t/2-(x_1-y_1)^2/2t}\\\label{j:estimates}
\approx&\frac{x_1y_1}{t+x_1y_1}\frac{(1-x_1)(1-y_1)}{t+(1-x_1)(1-y_1)}\frac{1+t^{5/2}}{\sqrt t}e^{-\pi^2t/2-(x_1-y_1)^2/2t}.
\end{align}
Combining (\ref{j:estimates}) with the formula (\ref{kBdens}), we obtain
\begin{align*}
G_{S_{a,1}}^{(\mu)}(x,y)\approx&\,x_1y_1(1-x_1)(1-y_1)(x_n-a)(y_n-a)x_n^{2\mu-1}y_n^{-2}\\
&\hspace{-14mm}\times\int_0^\infty\frac{1}{t+x_1y_1}\frac{1}{t+(1-x_1)(1-y_1)}\frac{1+t^{5/2}}{t^{n/2}}\frac{e^{-\pi^2t/2-\left|x-y\right|^2/2t}}{t+(x_n-a)(y_n-a)}\left(\frac{1}{t+x_ny_n}\right)^{\mu-\frac{1}{2}}dt.
\end{align*}
Next, we apply Lemma \ref{intlemma} with $\alpha=\frac52$, $\beta=\frac{n-2}{2}$, $b=\left|x-y\right|^2$, $k=4$ $a_1=x_1y_1$, $\gamma_1=1$, $a_2=(1-x_1)(1-y_1)$, $\gamma_2=1$, $a_3=(x_n-a)(y_n-a)$, $\gamma_3=1$, $a_4=x_ny_n$, $\gamma_4=\mu-\frac12$ and get
\begin{align*}
G_{S_{a,1}}^{(\mu)}(x,y)\stackrel{\mu,n}{\approx}
 \frac{x_n^{\mu-\frac12}}{y_n^{\mu+\frac32}}\frac{\delta_1(x_1)\delta_1(y_1)\(1+|x-y|\)^{\mu+\frac n2+\frac 72}e^{-\pi\left|x-y\right|^2}}{|x-y|^{n-2}\(\left|x-y\right|+\cosh\rho_a\)\(\left|x-y\right|+\cosh\rho\)^{\mu-\frac{1}{2}}}\,w(x,y),
\end{align*}
where 
\begin{align}\label{w}
w&(x,y)=\\\nonumber
&=\frac{1}{x_1y_1+x_1y_1\left|x-y\right|+\left|x-y\right|^2}\frac{1}{(1-x_1)(1-y_1)+(1-x_1)(1-y_1)\left|x-y\right|+\left|x-y\right|^2}.
\end{align}
To complete the proof we need to show that 
\begin{equation}\label{w:estimates}
w(x,y)\approx \frac{\(\delta_1(x_1)\delta_1(y_1)\)\wedge |x-y|^2}{\delta_1(x_1)\delta_1(y_1)|x-y|^2\(1+|x-y|^2\)}.
\end{equation}
For $|x_1-y_1|\geq\frac12$ we get $w(x,y)\approx\frac1{|x-y|^4}$, which is equivalent to   (\ref{w:estimates}) in this case. On the other hand, for $|x-y|<\frac12$ we  have $x_1y_1\approx1$ or $(1-x_1)(1-y_1)\approx1$ and consequently
\begin{align*}
w(x,y)&\approx
\frac{1}{x_1y_1+\left|x-y\right|^2}\frac{1}{(1-x_1)(1-y_1)+\left|x-y\right|^2}\\
&\approx
\frac{1}{[(1-x_1)(1-y_1)]\wedge[x_1y_1]+\left|x-y\right|^2}\\
&\approx
\frac{1}{\delta_1(x_1)\delta_1(y_1)+\left|x-y\right|^2},
\end{align*}
as required. 
\end{proof}

The Poisson kernel of smooth and bounded domains may be obtained as a derivative of the Green function with respect to the normal vector. Since the set  $S_{a,b}$ is neither bounded  nor smooth, we derive its Poisson kernel separately.

\begin{theorem}\label{estP}For $x\in S_{a,b}$, $y\in \partial S_{a,b}$ we have
\begin{align*}P_{S_{a,b}}^{(\mu)}(x,y)&\stackrel{\mu,n}{\approx}\(\dfrac{x_n}{y_n}\)^{\mu-1/2}\frac{e^{-\frac\pi b\left|x-y\right|}\(1+\frac1b|x-y|\)^{\mu+(n+3)/2}}{|x-y|^n\(\frac1b\left|x-y\right|+\cosh\rho\)^{\mu-1/2}}\\[10pt]
&\hspace{1cm}\times\left\{\begin{array}{ll}\dfrac{\delta_b(x_1)}{\(\frac1b\left|x-y\right|+\cosh\rho_a\)},&y_1\in
\{0,b\},\\[15pt]
(x_n-y_n)\dfrac{\[\delta_b(x_1)\delta_b(y_1)\]\wedge |x-y|^2}{|x-y|^{2}},&y_n=a.\end{array}\right.
\end{align*}
\end{theorem}
\begin{proof} In view of scaling property and  Theorem \ref{XtoY} we need only  to investigate the density function of $Y\(\tau^Y_{S_{a,1}}\)$.   Let $\tau^B_{(0,1)}$ be the first exit time from $(0,1)$ by the Brownian motion $B_1(t)$ and $\tau^R_{(a,\infty)}$ be the first exit time from $(a,\infty)$ by the Bessel process $R^{(-\mu)}(t)$. Observe that $$\tau^Y_{S_{a,1}}=\tau^B_{(0,1)}\wedge\tau^R_{(a,\infty)}.$$
Furthermore, let us divide the boundary $\partial S_{a,1}$ of $S_{a,1}$ into two parts: $\partial_1 S_{a,1}=\{0,1\}\times\R^{n-2}\times(a,\infty)$ and $\partial_2S_{a,1}=(0,1)\times\R^{n-2}\times\{a\}$. For any Borel set $A\subset \partial_1 S_{a,1}$ we have 
\begin{align*}
\p^x&\(Y\(\tau^Y_{S_{a,1}}\)\in A\)=\ \p^x\(\(B_1\(\tau^B_{(0,1)}\),...,B_{n-1}\(\tau^B_{(0,1)}\),R^{(-\mu)}\(\tau^B_{(0,1)}\)\)\in A,\tau^B_{(0,1)}\leq\tau^R_{(a,\infty)}\).
\end{align*}
Since $\tau^B_{(0,1)}$ is independent of the rest of the above-appearing processes and variables, we may write
$$\p^x\(Y\(\tau^Y_{S_{a,1}}\)\in A\)=\int_A\int_0^\infty
\gamma(t;x_1,y_1)\frac{\exp\(-\frac1{2t}\sum_{k=2}^{n-1}(x_k-y_k)^2\)}{(2\pi
t)^{(n-2)/2}}g^{(\mu)}_a(t;x_n,y_n)dt\,dy,$$
where $\gamma(t;x_1,y_1)=\p^x\(\beta_1\(\tau^\beta_{(0,1)}\)=y_1, \tau^\beta_{(0,1)}\in dt \)/dt$. Consequently, the inner integral represents the Poisson kernel $P_{S_{a,1}}^{(\mu)}(x,y)$.  Using the following estimates of  the function $\gamma(t;x_1,y_1)$  (see \cite{PSZ}, Thm. 5.3)
$$\gamma(t;x_1,y_1)\,\approx\,
x_1(1-x_1)\frac{1+t^{5/2}}{(t+1-|x_1-y_1|)t^{3/2}}\exp\(-\frac{|x_1-y_1|^2}{2t}-\frac12\pi^2t\),$$
where $x_1,y_1\in (0,1), t>0$, and the formula (\ref{kBdens}), we obtain
\begin{align*}
P_{S_{a,1}}^{(\mu)}(x,y)\stackrel{\mu,n}{\approx}& x_1(1-x_1)(x_n-a)(y_n-a)x_n^{2\mu-1}\\
&\times\int_0^\infty
\frac{1+t^{5/2}}{(t+1-|x_1-y_1|)t^{(n+2)/2}}\frac{\exp\(-\frac{\left|x-y\right|^2}{2t}-\frac12\pi^2t\)}{t+(x_n-a)(y_n-a)}\left(\frac{1}{t+x_ny_n}\right)^{\mu-\frac{1}{2}}dt,
\end{align*}
To deal with this integral, we  apply Lemma \ref{intlemma} with $\alpha=\frac52$, $\beta=\frac{n}{2}$, $b=\left|x-y\right|^2$, $k=3$ $a_1=1-|x_1-y_1|$, $\gamma_1=1$, $a_2=(x_n-a)(y_n-a)$, $\gamma_2=1$, $a_3=x_ny_n$, $\gamma_3=\mu-\frac12$ and get
\begin{align*}
P_{S_{a,1}}^{(\mu)}(x,y)\stackrel{\mu,n}{\approx}&\, \frac{x_1(1-x_1)e^{-\pi\left|x-y\right|}\(1+|x-y|\)^{\mu+(n+7)/2}}{|x-y|^n\(1+\left|x-y\right|+\frac{\left|x-y\right|^2}{(x_n-a)(y_n-a)}\)\(1+\left|x-y\right|+\frac{\left|x-y\right|^2}{x_ny_n}\)^{\mu-1/2}}\\
&\times \(\frac{x_n}{y_n}\)^{\mu-1/2}\frac{1}{1-|x_1-y_1|+(1-|x_1-y_1|)\left|x-y\right|+\left|x-y\right|^2}.
\end{align*}
 Using $1+|x-y|-|x_1-y_1|\approx 1+|x-y|$, we estimate  the denominator of the last fraction as follows
\begin{align*}
1-|x_1-y_1|+(1-&|x_1-y_1|)\left|x-y\right|+\left|x-y\right|^2\\
=&(1-|x_1-y_1|+|x-y|)-|x_1-y_1|\left|x-y\right|+\left|x-y\right|^2\\
\approx&1+|x-y|(1-|x_1-y_1|+\left|x-y\right|)\\
\approx&1+|x-y|^2.
\end{align*}
Eventually we arrive at
\begin{align*}
P_{S_{a,1}}^{(\mu)}(x,y)\stackrel{\mu,n}{\approx} \(\frac{x_n}{y_n}\)^{\mu-1/2}\frac{\delta_1(x_1)e^{-\pi\left|x-y\right|}\(1+|x-y|\)^{\mu+(n+3)/2}}{|x-y|^n\(\left|x-y\right|+\cosh\rho_a\)\(\left|x-y\right|+\cosh\rho\)^{\mu-1/2}},
\end{align*}
which completes estimates of the Poisson kernel on the first part of the boundary. Assume now  $B\subset\partial_2S_{a,1}$. Note that $y_n=a$ for $y\in \partial_2S_{a,1}$. Similarly as in previous case we get
\begin{align*}
\p^x\(Y\(\tau^Y_{S_{a,1}}\)\in A\)=&\ \p^x\(\(\beta_1\(\tau^R_{(a,\infty)}\),...,\beta_{n-1}\(\tau^R_{(a,\infty)}\),a\)\in A,\tau^R_{(a,\infty)} <\tau^\beta_{(0,1)}\)\\
=&\ \int_B\int_0^\infty
j(t;x_1,y_1)\frac{\exp\(-\frac1{2t}\sum_{k=2}^{n-1}(x_k-y_k)^2\)}{(2\pi
t)^{(n-2)/2}}q^{(-\mu)}_{a}(t;x_n)dt\,dy,
\end{align*}
where  $q^{(-\mu)}_{a}(t;x_n)=\p^{x_n}\(\tau^R_{(a,\infty)}\in dt\)/dt$. 
Hence, by (\ref{j:estimates}) and (\ref{kBtime}), we obtain
\begin{align*}
P_{S_{a,1}}^{(\mu)}(x,y)\stackrel{\mu,n}{\approx}&x_1y_1(1-x_1)(1-y_1)(x_n-y_n)x_n^{2\mu-1}\\
&\times\int_0^\infty
\frac{1+t^{5/2}}{t+x_1y_1}\frac{e^{-\pi^2t/2-|x-y|^2/2t}}{t+(1-x_1)(1-y_1)}\frac{1}{
t^{(n+2)/2}}\frac{1  }{\(t+ y_nx_n\)^{\mu-1/2}}dt.
\end{align*}
Next, we apply Lemma \ref{intlemma} with $\alpha=\frac52$, $\beta=\frac{n}{2}$, $b=\left|x-y\right|^2$, $k=3$, $a_1=x_1y_1$, $\gamma_1=1$, $a_2=(1-x_1)(1-y_1)$, $\gamma_2=1$, $a_3=y_nx_n$, $\gamma_3=\mu-\frac12$ and get
\begin{align*}
P_{S_{a,1}}^{(\mu)}(x,y)\stackrel{\mu,n}{\approx}&\frac{x_1y_1(1-x_1)(1-y_1)(x_n-y_n)x^{2\mu-1}}{\(y_nx_n+y_nx_n|x-y|+|x-y|^2\)^{\mu-1/2}e^{-\pi|x-y|}}\frac{\(1+|x-y|\)^{\mu+(n+7)/2}}{|x-y|^n}w(x,y)\\
\approx&\(\frac{x_n}{y_n}\)^{\mu-1/2}\frac{\delta_1(x_1)\delta_1(y_1)(x_n-y_n)e^{-\pi|x-y|}}{|x-y|^n}\frac{\(1+|x-y|\)^{\mu+(n+7)/2}}{\(|x-y|+\cosh\rho\)^{\mu-1/2}}w(x,y),
\end{align*}where $w(x,y)$ is given by (\ref{w}). Usage of the estimate (\ref{w:estimates}) ends the proof.
\end{proof}
Manipulating with parameters $a$ and $b$ in Theorems \ref{estG} and \ref{estP}, we obtain some further results. Calculating limits as $a\rightarrow0$ and using monotone convergence theorem we get the below-given corollary. It generalizes  estimates from \cite{S}, where  HBM without drift was considered. 
\begin{corollary}
For $x,y\in S_{0,b}$ we have
$$
G^{(\mu)}_{S_{0,b}}(x,y)\stackrel{\mu,n}{\approx}\frac{x_n^{\mu-1/2}}{y_n^{\mu+3/2}}\frac{e^{-\frac{\pi}b\left|x-y\right|}}{\left|x-y\right|^{n}}\(\delta_b(x_1)\delta_b(y_1)\wedge\left|x-y\right|^2\)\frac{\(1+\frac1b\left|x-y\right|\)^{n/2+\mu+3/2}}{\(\frac1b\left|x-y\right|+\cosh\rho\)^{\mu+1/2}},
$$
and for  $x\in S_{0,b}$, $y\in \partial S_{0,b}$ we have
\begin{align*}
P_{S_{0,b}}^{(\mu)}&(x,y)\\
&\stackrel{\mu,n}{\approx}\left\{\begin{array}{ll} \(\dfrac{x_n}{y_n}\)^{\mu-1/2}\dfrac{\delta_b(x_1)e^{-\frac\pi b\left|x-y\right|}\(1+\frac1b|x-y|\)^{\mu+(n+3)/2}}{|x-y|^n\(\frac1b\left|x-y\right|+\cosh\rho\)^{\mu+1/2}},&y_1\in
\{0,b\},\\[15pt]
x_n^{2\mu}e^{-\frac\pi b\left|x-y\right|}\Big(\delta_b(x_1)\delta_b(y_1)\wedge |x-y|^2\Big)\dfrac{\(1+\frac1b|x-y|\)^{\mu+(n+3)/2}}{|x-y|^{2\mu+n+1}},&y_n=a.\end{array}\right.
\end{align*}
\end{corollary}
Taking additionally limits as $b\rightarrow\infty$ we obtain estimates provided in  \cite{MS}. The next corollary concerns  the mostly studied   subset of $\H^n$ in context of HBM i.e. $D_a=\{x\in\H^n:x_n>a\}$, $a>0$. It follows from Theorems \ref{estG} and \ref{estP} by replacing $x_1$ and $y_1$ by $x_1+\frac b2$ and  $y_1+\frac b2$, respectively, and   taking limits as $b$ tends to infinity.  In fact,  the Poisson kernel $P_{D_a}^{(\mu)}(x,y)$ was estimated in \cite{BMR}, and  estimates of the $\lambda$-Green function for the process without drift  (which, by Theorem \ref{thm:allis0}, are equivalent to estimates of the Green function for the process with suitable drift) are the main results of   \cite{BBM}.
\begin{corollary}
For $\mu>0$ we have
$$\begin{array}{rll}
G^{(\mu)}_{D_a}(x,y)\stackrel{\mu,n}{\approx}&\dfrac{x_n^{\mu-1/2}}{y_n^{\mu+3/2}}\dfrac{1}{\left|x-y\right|^{n-2}\cosh\rho_a\(\cosh\rho\)^{\mu-1/2}},&x,y\in D_a,\\[15pt] 
P_{D_a}^{(\mu)}(x,y)\stackrel{\mu,n}{\approx}&\(\dfrac{x_n}{y_n}\)^{\mu-1/2}\dfrac{x_n-y_n}{|x-y|^n\(\cosh\rho\)^{\mu-1/2}},&x\in D_a, y\in \partial D_a.
\end{array}$$
\end{corollary}

\section{Appendix}
In this section we present a technical lemma  which is used to estimate integrals appearing in Section 4.
\begin{lemma}\label{intlemma}
Fix $\alpha\geq0$, $\beta\geq\frac12$, $k\in\{0,1,2,...\}$ and $\gamma_i\geq0$, $i\in\{1,...,k\}$. There exists a constant $c=c(\alpha,
\beta,
\gamma_1,...,\gamma_k)$ such that for $a_i>0$, $i\in\{1,...,k\}$, and $b>0$ we have
$$\int_0^\infty
\frac{(1+t)^{\alpha}}{t^{\beta+1}}\frac{\exp\(-\frac{b^2}{2t}-\frac12\pi^2t\)}{\prod_{i=1}^k\(a_i+t\)^{\gamma_i}}dt\stackrel{c}{\approx}
\frac{e^{-b\pi}}{b^{2\beta}}\frac{1+b^{\alpha+\beta-1/2+\sum_{i=1}^k\gamma_i}}{\prod_{i=1}^k\(a_i+a_ib+b^2\)^{\gamma_i}}.$$
Additionally, the estimates stay valid also if there is one index $i\in\{1,...,k\}$  such that $\gamma_i$ is negative but greater than $-\frac12$.
\end{lemma}
\begin{proof}Throughout this proof only, every letter  $c$ appearing over the sign $\approx$ represents a constant depending on all of parameters: $\alpha$, $\beta$, $\gamma_1,...,\gamma_k$.  Substituting $t=\frac{bu}\pi$ in the integral from the thesis we get
$$
\frac{e^{-b\pi}}{b^{\beta}}\pi^{\beta-\alpha+\sum_{i=1}^k\gamma_i}\int_0^\infty
\frac{(\pi+ub)^{\alpha}}{u^{\beta+1}}\frac{\exp\(-\frac12b\pi\(\sqrt u - \frac1{\sqrt
u}\)^2\)}{\prod_{i=1}^k\(a_i\pi+bu\)^{\gamma_i}}du=:\frac{e^{-b\pi}}{b^{\beta}}\mathcal I.
$$
It is  enough to show 
\begin{equation}\label{I}
\mathcal I\stackrel{c}{\approx}\frac{1+b^{\alpha+\beta-1/2+\sum_{i=1}^k\gamma_i}}{b^{\beta}\prod_{i=1}^k\(a_i+a_ib+b^2\)^{\gamma_i}}.
\end{equation}
We start with the substitution  $\sqrt u-\frac{1}{\sqrt u}=s\sqrt{\frac2b}$. Note that
$$u=\(\sqrt{1+\frac{s^2}{2b}}+\frac{s}{\sqrt {2b}}\)^2{\approx}\left\{\begin{array}{ll}1+\dfrac{s^2}{b},\ \ \
&s>0\  (\Leftrightarrow u>1),\\[12pt]\dfrac1{1+\frac{s^2}{b}},&s\leq0\ (\Leftrightarrow u\leq1),\\
\end{array}\right.$$
and $$\frac{du}{u}=\frac{2ds}{\sqrt{s^2+2b}}{\approx}\frac{ds}{\sqrt{s^2+b}}.$$
Consequently we obtain
\begin{eqnarray*}
\mathcal I&\stackrel{c}{\approx}& b^{\beta}\int_0^\infty \frac{(1+b+
s^2)^{\alpha}}{\(b+s^2\)^{\beta+1/2}}\frac{e^{-s^2}}{\prod_{i=1}^k\(a_i+b+s\)^{\gamma_i}}ds+\\
&&+\frac{1}{b^{\beta}}\int_{-\infty}^0 \frac{\(s^2+b\)^{\beta-1/2}\(1+\frac{b^2
}{s^2+b}\)^{\alpha}e^{-s^2}}{\prod_{i=1}^k\(a_i+\frac{b^2
}{s^2+b}\)^{\gamma_i}}ds\\
&=:&b^{\beta} \mathcal I_1+\frac{1}{b^{\beta}} \mathcal I_2.
\end{eqnarray*}
For $b>1$ we have
\begin{eqnarray*}
\mathcal I_1&=&\frac{(1+b)^{\alpha}}{b^{\beta+1/2}\prod_{i=1}^k(a_i+b)^{\gamma_i}}\int_0^\infty \frac{(1+
\frac{s^2}{1+b})^{\alpha-1}}{\(1+\frac{s^2}{b}\)^{\beta+1/2}}\frac{e^{-s^2}}{\prod_{i=1}^k\(1+\frac{s^2}{a_i+b}\)^{\gamma_i}}ds\\
&\stackrel{c}{\approx}&\frac{b^{-\beta+\alpha-1/2}}{\prod_{i=1}^k(a_i+b)^{\gamma_i}},\\
\mathcal I_2&<&\int_{-\infty}^0
\frac{\(b(s^2+1)\)^{\beta-1/2}(1+b)^{\alpha}e^{-s^2}}{\prod_{i=1}^k\(\frac{a_i+b
}{s^2+1}\)^{\gamma_i}}ds\stackrel{c}{\approx}\frac{b^{\beta+\alpha-1/2}}{\prod_{i=1}^k(a_i+b)^{\gamma_i}}.
\end{eqnarray*}
Hence
\begin{equation}\label{I>1} \mathcal I\stackrel{c}{\approx}b^\beta \mathcal I_1\stackrel{c}{\approx}
\frac{b^{\alpha-1/2}}{\prod_{i=1}^k(a_i+b)^{\gamma_i}},\ \ \ \ \ \ b>1,\end{equation}
which is equivalent to (\ref{I}).
Let assume now $b\leq1$. We have
$$\mathcal I_2\stackrel{c}{\approx}\int_{-\infty}^0
\frac{(s^2+b)^{\beta-1/2}e^{-s^2}}{\prod_{i=1}^k\(a_i+\frac{b^2
}{s^2+b}\)^{\gamma_i}}ds.$$
We are going now to use inequalities $s^2<s^2+b<s^2+1$ and $\(a_i+b^2\)\frac{1}{1+s^2}<a_i+\frac{b^2
}{s^2+b}<\(a_i+b^2\)\(1+\frac1{s^2}\)$, $1\leq i\leq k$. Note that, as long as all of $\gamma_i$ are nonnegative, replacing $a_i+\frac{b^2}{s^2+b}$ by $\(a_i+b^2\)\(1+\frac1{s^2}\)$ does not change integrability of the above integral. However, it works also if one of $\gamma_i$ belongs to the interval $\(-\frac12,0\)$. Thus
$$\mathcal I_2\stackrel{c}{\approx}\frac{1}{\prod_{i=1}^k\(a_i+b^2\)^{\gamma_i}},\ \ \ \ b\leq1.$$
Furthermore, we estimate $\mathcal I_1$ as follows
\begin{eqnarray*}
\mathcal I_1&\stackrel{c}{\approx}&\frac{1}{b^{\beta+1/2}\prod_{i=1}^k\(a_i+b^2\)^{\gamma_i}}\int_0^\infty
\frac{(1+s^2)^{\alpha}}{\(1+\frac{s^2}{b}\)^{\beta+1/2}}\frac{e^{-s^2}}{\prod_{i=1}^k\(1+\frac {s^2}{a_i+b^2
}\)^{\gamma_i}}ds\\
&<&\frac{1}{b^{\beta+1/2}\prod_{i=1}^k\(a_i+b^2\)^{\gamma_i}}\int_0^\infty (1+s^2)^{\alpha}e^{-s^2}ds\\
&\stackrel{c}{\approx}&\frac{1}{b^{\beta+1/2}\prod_{i=1}^k\(a_i+b^2\)^{\gamma_i}}.
\end{eqnarray*}
Finally we get
$$\mathcal I\stackrel{c}{\approx}\frac1{b^\beta}\,\mathcal
I_2\stackrel{c}{\approx}\frac{1}{b^\beta\prod_{i=1}^k\(a_i+b^2\)^{\gamma_i}},\ \ \ \ b\leq1.$$
This coincides with  (\ref{I}) and the proof is complete.
\end{proof}

\end{document}